\documentclass[11pt]{amsart}
 
\usepackage{amsmath,amsthm,amssymb}
\usepackage{graphicx,psfrag,epsfig}

\begin{document}
\title{A geometric proof of the existence of Whitney stratifications}
\author{V.~Yu.~Kaloshin \thanks{Department of Mathematics,
Princeton University, Princeton NJ 08544-1000}}
\thanks{The author is partially supported by the Sloan Dissertation
Fellowship and the American Institute of Mathematics Five-Year Fellowship}
\email{kaloshin@math.princeton.edu}

\maketitle
\markboth{Whitney's stratification}{V. Yu. Kaloshin}       
 
\theoremstyle{plain}
\newtheorem{Thm}{Theorem} 
\newtheorem{Def}{Definition}
\newtheorem{Lm}{Lemma}                             
\newtheorem{Rem}{Remark}

\def\bdef{\begin{Def}}
\def\endef{\end{Def}}
\def\bthm{\begin{Thm}}
\def\ethm{\end{Thm}} 
\def\blm{\begin{Lm}}
\def\elm{\end{Lm}}
\def\brm{\begin{Rem}}
\def\erm{\end{Rem}}
\def\beq{\begin{eqnarray}}
\def\eneq{\end{eqnarray}}
\def\Cal{\mathcal}
\def\Bbb{\mathbb}
\def\gm{\gamma}
\def\dt{\delta}
\def\eps{\epsilon}
\def\lb{\lambda}
\def\Bbb{\mathbb}
\def\R{\Bbb R}
\def\P{\Cal P}
\def\C{\Bbb C}
\def\Z{\Bbb Z}

\section{Introduction}

A stratification of a set, e.g. an analytic variety, is, 
roughly, a partition of it into manifolds so that these manifolds
fit together ``regularly''. Stratification theory was originated 
by Thom and Whitney for algebraic and analytic sets.
It was one of the key ingredients in Mather's proof of
the topological stability theorem {\cite{Ma}} 
(see {\cite{GM}} and {\cite{PW}} for the history and further 
applications of stratification theory). 

In this paper, given a partition of a singular set (which we know
always exists), we prove that there is a ``regular'' partition. Our 
proof is based on a remark that if there are two parts of the partition 
$V$ and $W$ of different dimension and $V \subset \overline{W}$,  then 
{\it irregularity} of the partition at a point $x$ in $V$ corresponds to
the existence of {\it nonunique} limits of tangent planes $T_yW$ as $y$ 
approaches $x$. 

Consider either the category of (semi)analytic (or (semi)algebraic) 
sets. Call a subset $V \subset \R^m$ (or $\C^m$) a {\it{semivariety}}
if locally at each point $x\in \R^m$ (or $\C^m$) it is a finite
union of subsets defined by equations and inequalities
\beq
f_1= \dots =f_k=0 \quad 
\begin{cases}
g_1\neq 0, \dots, g_l \neq 0\ \ \  (\textup{complex case}),\\
g_1>0, \dots, g_l> 0\ \ \  (\textup{real case}),
\end{cases}
\eneq
where $f_i$'s and $g_j$'s are real (or complex) analytic 
(or algebraic) depending on the case under consideration.

In the real algebraic case semivarieties are usually called 
semialgebraic sets; in the complex algebraic case they are called
constructible, and in either analytic case they are called 
semianalytic sets. Semivarieties are closed under
Boolean operations. 

\bdef  \label{regular} (Whitney)
Let $V_i,V_j$ be disjoint manifolds in $\R^m$ (or $\C^m$),
$\dim V_j>\dim V_i$, and let $x \in V_i \cap \overline{V_j}$. 
A triple $(V_j,V_i,x)$ is called $a$(resp. $b$)- regular if 

$A$) when a sequence $\{y_n\} \subset V_j$ tends to $x$
and $T_{y_n} V_j$ tends in the Grassmanian bundle to 
a subspace $\tau_x$ of $\R^m$ (or $\C^m$), then 
$T_x V_i\subset \tau_x$;

$B$) when sequences $\{y_n\} \subset V_j$ and 
$\{x_n\} \subset V_i$ each tends to $x$, the unit vector 
$(x_n-y_n)/|x_n-y_n|$ tends to a vector $v$, and 
$T_{y_n}V_j$ tends to $\tau_x$, then 
$v\in \tau_x$\footnote{This way of defining $b$-regularity is 
due to Mather \cite{Ma}. Whitney's definition \cite{Wh}
is equivalent to this one provided of $a$-regularity}.

$V_j$ is called $a$(resp. $b$)- regular over 
$V_i$ if each triple $(V_j,V_i,x)$ is $a$(resp. $b$)- regular.
\endef

\bdef (Whitney) Let $V$ be a semivariety in $\R^m$ (or $\C^m$). 
A disjoint decomposition 
\beq\label{decomp}
V=\bigsqcup_{i \in I} V_i, \quad V_i\cup V_j =\emptyset \ \ \ 
\textup{for} \ \ \ i\neq j  
\eneq
into smooth semivarieties $\{V_i\}_{i \in I}$,
called strata, is called an $a$(resp. $b$)-regular stratification  
if

1. each point has a neighborhood intersecting only finitely 
many strata;

2. the frontier $\overline{V_j}\setminus V_j$ of each stratum $V_j$ 
is a union of other strata $\bigsqcup_{i \in J(i)} V_i$;

3. any triple $(V_j,V_i,x)$ such that $x \in V_i \subset
\overline{V_j}$ is $a$(resp. $b$)-regular.
\endef

\bthm \label{main} \cite{Wh},\cite{Th},{\cite{Lo}} For any semivariety 
$V$ in $\R^m$ (or $\C^m$) there is an $a$(resp. $b$)-regular
stratification.
\ethm
The existence of stratifications in the complex analytic case
was proved by Whitney {\cite{Wh}}. Later Thom published a sketch 
of a proof {\cite{Th}}. Then Lojasiewicz {\cite{Lo}} extended these
results to the semianalytic case. The most illuminating proof is due 
to Wall {\cite{Wa}}, where based on Milnor's curve selection lemma 
{\cite{Mi}} he simplifies the above proofs. Hironaka {\cite{Hi}} gave 
an elegant proof using his resolution of singularities, but 
it requires background in algebraic geometry. We give 
a geometric proof based on Milnor's curve selection lemma \cite{Mi},
\cite{Wa}, Rolle's lemma, and a transversality 
theorem. The rest of the paper is devoted to this proof.

{\it{Proof of theorem \ref{main}:}}\ \ A semivariety $V$ has
well-defined dimension, say $d \leq m$. Denote by $V_{reg}$ the set of
points, where  $V$ is locally a real (or complex) analytic submanifold 
of $\R^m$ (or $\C^m$) of dimension $d$. $V_{reg}$ is a semivariety, 
moreover, $V_{sing}=V \setminus V_{reg}$ is a semivariety of positive
codimension in $V$, i.e. $\dim V_{sing}< \dim V$. In the analytic case
all these results may be found in Lojasiewicz {\cite{Lo}}; in the
algebraic case they are not difficult (see e.g. {\cite{Mi}}).

{\it{Step 1.}} There is a filtration of $V$ by semivarieties
\beq
V^0 \subset V^1 \subset \dots \subset V^d=V,
\eneq
where for each $k=1, \dots, d$ the set $V^k \setminus V^{k-1}$ 
is a manifold of dimension $k$. This follows from the Lojasiewicz 
result. Indeed, consider $V_{sing} \subset V$, then 
$V\setminus V_{sing}$ is a manifold of dimension $d$ and 
$\dim V_{sing}<d$. Inductive application of these arguments
completes the proof.

A {\it{refinement}} of a decomposition $V=\bigsqcup_{i\in I}V_i$
is a decomposition $V=\bigsqcup_{i'\in I'}V_{i'}$ such that 
any stratum $V_j$ of the first decomposition is a union of some 
strata of the second one, i.e. there is a set
$I'(j) \subset I'$ such that $V_j=\bigsqcup_{i' \in I'(j)}V_{i'}$.  
 
{\it{Step 2}}. Let $V\subset \R^m$ (or $\C^m$) be a manifold and 
$W\subset V$ be a semivariety. Denote by $Int_V(W)$ the set of 
interior points of $W$ in $V$ w.r.t. the induced from $\R^m$ 
(resp. $\C^m$) topology. Let $V_i$ and $V_j$ be a pair of distinct
strata. For each point $x \in V_i\cap \overline{V_j}$ denote by 
$V_j^{con,x}$ a local connected component  of $V_j$ at $x$, i.e. 
a connected component of intersection of $V_j$ with a ball centered 
at $x$ and call it {\it essential} if the closure of $V_j^{con,x}$
has $x$ is in the interior, 
$x\in Int_{V_i}(V_i\cap \overline{V^{con,x}_j})$. Denote by 
$V_j^{ess,x}$ the union of all local essential components of $V_j$.
Lojasiewicz {\cite{Lo}} showed that $V_j$ has only a finitely many 
local connected components. 
 
\bthm \label{sing} For any two disjoint strata $V_j$ and 
$V_i$ the set of points 
\beq \nonumber
Sing_{a(\textup{resp.} b)}(V_j,V_i)=
\{ x \in V_i\cap\overline{V_j}:\ (V_j^{ess,x},V_i,x) \ 
\textup{is not}\ a (\textup{resp.}\ b)-\textup{regular}\},
\eneq
is a semivariety in $V_i$ and 
$\dim  Sing_{a(\textup{resp.} b)}(V_j,V_i)<\dim V_i$.
\ethm
Let us show that this theorem is sufficient to prove Theorem
\ref{main}.  Consider a decomposition $V=\bigsqcup_{i\in I}V_i$ and
split the strata into two groups: the first group consists of strata 
of dimension at least $k$ and the second group is of the rest. 
Suppose that each stratum from the first group is $a$(resp. $b$)-regular
over each stratum from the second group. Then by definition of 
$a$(resp. $b$)-regularity any refinement of a stratum from the second 
group preserves this $a$(resp. $b$)-regularity. 

Now apply this refinement inductively. Consider strata in 
$V^d\setminus V^{d-1}$ of dimension $d$. Using Theorem \ref{sing} and
the result of Lojasiewicz \cite{Lo} that a frontier of a semivariety 
has dimension less than a semivariety itself, 
refine $V^{d-1}$ so that each $d$-dimensional stratum is 
$a(\textup{resp.}\ b)$-regular over each stratum in $V^{d-1}$. 
The above remark shows that any further refinement of the 
strata in $V^{d-1}$ preserves the $a(\textup{resp.}\ b)$-regularity
of strata from $V^d\setminus V^{d-1}$ over it. This reduces
the problem of the existence of stratification for $d$-dimensional
semivarieties to the same problem for $(d-1)$-dimensional semivarieties. 
Induction on dimension completes the proof of Theorem \ref{main}.

Our proof is based on the observation that if $V_i\subset \overline{V_j}$
are a pair of strata $a$(resp. $b$)-regularity of $V_j$ over $V_i$ at
$x$ in $V_i$ is closely related to whether the limit of tangent planes 
$T_yV_j$ is unique or not as $y$ from $V_j$ tends to $x$.
The rest of the paper is devoted to the proof of Theorem \ref{sing} 
which consists of two steps. In section \ref{key} we relate 
$a$(resp. $b$)-regularity with (non)uniqueness of limits of tangent
planes $T_yV_j$, then based on it and Rolle's lemma in section 
\ref{reduction} we prove Theorem \ref{sing}.

\subsection{The key definitions} \label{key}
Let $V_i$ and $V_j$ be a pair of distinct strata.  Define
\beq 
\begin{aligned}\label{uniquelim}
Un_{a}(V_j,V_i)=\{ x \in V_i\cap \overline{V_j}&:
\textup{for any}\ V_j^{con,x},\ \textup{there exists}\ 
\tau_x\subset T_x\R^m\\ (\textup{resp.}\ T_x\C^m)\ 
\textup{such that for any}\ & \{y_n\} \subset V_j^{con,x}\ 
\textup{tending to}\  x, \ \  T_{y_n}V_j \to \tau_x \},
\end{aligned}
\eneq

Since $a$(resp. $b$)-regularity is a local property, w.l.o.g. we can 
assume that locally $V_i$ is an $s$-plane with 
a basis of unit vectors $e_1,\dots ,e_s$.  Using an idea of Kuo 
{\cite{Ku}} (see also {\cite{Wa}}) we define {\it a Kuo map} 
$\P^{a(resp.\ b)}: V_j \to \R$ which measures non 
$a$(resp. $b$)-regularity in terms of an angle between 
a vector or a plane and the tangent plane to $V_j$. 
Denote by $\pi_i:\R^m \to V_i^\perp$ (resp. $\pi_i^\perp:\R^m\to V_i$)
the orthogonal projection along $V_i$ (resp. $V_i^\perp$) onto 
the complement $V_i^\perp$ (resp. $V_i$) with $x$ being the origin 
of $\R^m$ and by 
\beq 
\begin{aligned}\label{kuofunc}
\pi_j:V_j\times \R^m \to \R^m, \ \pi_{j,t}: V_j \to \R^m \ 
\textup{for}\ t=1,\dots,s+1\ \textup{defined by}\\ 
\pi_j(y,v)=\pi_{T_yV_j}(v), \  \pi_{j,t}(y)=\pi_j(y,e_t),\ 
\pi_{j,s+1}(y)=\pi_j(y,\pi_i(y)/|\pi_i(y)|),
\end{aligned}
\eneq
where $\pi_j(y,v)$ is the orthogonal projection of $v$ along the
tangent plane $T_yV_j$ to $V_j$ at $y$ naturally embedded into
$\R^m$. Define analytic functions 
$\P^{a(resp.\ b)}: V_j \to \R$ by 
$\P^a(y)=\sum_{t=1}^{s}|\pi_{j,t}(y)|^2$ (resp. 
$\P^b(y)= \sum_{t=1}^{s+1}|\pi_{j,t}(y)|^2$). By the definition 
the level sets of $\P^{a(resp.\ b)}$ are semivarieties. 

Notice that the first $s$ terms of the function $\P^a(y)$
measure the angle between $T_xV_i=V_i$ and $T_yV_j$ and the last
term measures the angle between the $V^\perp_i$- component
of $(y-x)/|y-x|$ and $T_yV_j$. Since any vector can be decomposed
into $V_i$ and $V_i^\perp$ components, this proves the following  

{\bf Fact 1.} {\it For any pair distinct strata $V_j$ and $V_i$ 
existence of a sequence $\{y_n\}\subset V_j$  tending to $x$ with
a nonzero limit of $\P^{a(resp.\ b)}(y_n)$ is equivalent to 
$a$(resp. $b$)-irregularity of $V_j$ over $V_i$ at $x$.}

\beq 
\begin{aligned}
Un_b(V_j,V_i)=\{ x \in Un_a(V_j,V_i):
\textup{for any}\ V_j^{con,x},\ \textup{there exists}\ 
\eps \in \R \\ \textup{such that for any}\ 
\{y_n\} \subset V_j^{con,x}\ \textup{tending to}\  x, 
\ \P^b(y_n) \to \eps \},
\end{aligned}
\eneq

\blm \label{one}
Let $V_i$ and $V_j$ be a pair of disjoint strata in $\R^m$ 
(or $\C^m$) with $V_i \cap \overline{V_j}\neq \emptyset$. 
Then $Sing_{a(\textup{resp.} b)}(V_j,V_i)$ and 
$Un_{a(\textup{resp.} b)}(V_j,V_i)$ are semivarieties and 
\beq\nonumber
Sing_{a}(V_j,V_i)\subset Sing_{b}(V_j,V_i),\quad 
Sing_{a(\textup{resp.} b)}(V_j,V_i) \subset 
V_i\setminus Un_{a(\textup{resp.} b)}(V_j,V_i).
\eneq
\elm
\brm
The new result here is that $Sing_{a(\textup{resp.} b)}(V_j,V_i)\subset 
V_i\setminus Un_{a(\textup{resp.} b)}(V_j,V_i)$. The other inclusion 
may be found in \cite{Wh}, \cite{Ma}, \cite{Lo}. 
\erm
{\it Proof:}\ \ \  Let's first prove that $Sing_a(V_j,V_i)$ 
is a semivariety. Consider
$V_i\times TV_j=\{(x,y,T_yV_j):\ x\in V_i,\ y \in V_j\}$.
It is a semivariety in an appropriate Grassmanian bundle
over $\R^m \times \R^m$ (resp. $\C^m \times \C^m$) and so 
is its closure. The condition $T_xV_i \not\subset\tau_x$
is semialgebraic and a projection of a semivariety is a semivariety. 
In the real (resp. complex) algebraic case it is called
the Tarski-Seidenberg Principle {\cite{Ja}} (resp. elimination
theory {\cite{Mu}}). In the real analytic case it depends on 
a generalization due to Lojasiewicz {\cite{Lo}} to
varieties analytic in some variables and algebraic in others.
In the complex analytic case, a proof may be found in {\cite{Wh}}.  
Similar arguments show $Sing_b(V_j,V_i)$ and 
$Un_{a(\textup{resp.} b)}(V_j,V_i)$ are semivarieties.

Now let's see that $Sing_a(V_j,V_i)\subset Sing_b(V_j,V_i)$.
For any sequence $\{y_n\} \subset V_j$ such that $T_{y_n}V_j$ 
has a limit $\tau_x$ as $y_n$ tends to  $x$ and any $v \in T_xV_i$ 
there is a sequence $\{x_n\}\subset V_i$ such that $x_n$ tends to 
$x$ slower than the sequence $\{y_n\}$, i.e. $|y_n-x|/|x_n-x|\to 0$ 
and the unit vectors $(x_n-y_n)/|x_n-y_n|$ tends to $v$ as 
$n\to \infty$\footnote{This was first noticed by J.Mather {\cite{Ma}}}. 
If $x\notin Sing_b(V_j,V_i)$, then $v$ belongs to $\tau_x$. 
Since any $v \in T_xV_i$ belongs to $\tau_x$, $T_xV_i$ also 
belongs to $\tau_x$.

To see that $Sing_a(V_j,V_i)\subset V_i\setminus Un_a(V_j,V_i)$,
suppose $x \in Sing_a(V_j,V_i) \cap Un_a(V_j,V_i)$. Fix an $a$-irregular 
essential local connected component $V^{con,x}_j$ of $V_j$ at $x$. 
There is a $\dim V_j$-plane $\tau_x$ such that for any sequence 
$\{y_n\}\subset V^{con,x}_j$ tending to $x$ we have 
$T_{y_n}V_j\to \tau_x$. Since $x\in Sing_a(V_j,V_i)$,
we have $T_xV_i\not\subset \tau_x$, i.e. there is a unit vector
$v\in T_xV_i$ which has a positive angle with $\tau_x$, i.e.
$<(v,\tau_x)=2\dt>0$. Denote by 
$C_{\dt,v}(x)=\{y\in \R^m:\ (\frac{y-x}{|y-x|},v)>1-\dt\}$
the $\dt$-cone around $v$ centered at $x$ and by $l_v(x)$ the 
ray starting at $x$ in the direction of $v$. The intersection
$V^{con,x}_j\cap C_{\dt,v}(x)=V^{con,x}_{j,\dt,v}$ is 
a semivariety and $l_v(x)$ is in its closure. By the Lojasiewicz
result $V^{con,x}_{j,\dt,v}$ consists of a finite number of
connected components. So one can choose a connected component
$W^{con,x}_{j,\dt,v}\subset V^{con,x}_{j,\dt,v}$ which contains
$l_v(x)$ in the closure. By Milnor's curve selection lemma
{\cite{Mi}}, {\cite{Wa}} there is an analytic curve
$\gm$ which belongs to $W^{con,x}_{j,\dt,v} \cup \{x\}$.
Since $\gm$ is analytic, it has a limiting tangent vector $w$ at $x$
which is by our construction should belong to $\tau_x$ and
$C_{\dt,v}(x)$. This is a contradiction with $<(v,\tau_x)=2\dt$. 

To see that $Sing_b(V_j,V_i)\subset V_i \setminus Un_b(V_j,V_i)$ 
it is sufficient to prove that $Sing_b(V_j,V_i)$ $\cap 
Un_a(V_j,V_i)\subset V_i \setminus Un_b(V_j,V_i)$. Let 
$x\in Sing_b(V_j,V_i)\cap Un_a(V_j,V_i)$ and $V_j^{con,x}$ be a 
$b$-irregular essential local connected component at $x$. Since 
$x\in Un_a(V_j,V_i)$, there is a unique limiting tangent plane 
$\tau_x=\lim T_{y_n}V_j$ independent of $\{y_n\}\subset V_j^{con,x}$ 
tending to $x$ and by the previous passage $x$ is $a$-regular, i.e. 
$V_i\subset \tau_x$. By Fact 1 and $b$-irregularity of $x$ there is 
a sequence $\{y_n\}\in V^{con,x}_j$ such that $|\P^b(y_n)|\to 2\dt\neq 0$. 
Let's prove existence of a sequence $\{y'_n\}\in V^{con,x}_j$ such 
that $|\P^b(y_n)|\to \eps<\dt$ which shows that $x\notin Un_b(V_j,V_i)$.

For each $\tilde x \in V_i$ close to $x$ consider the ``level''
set $V_j^{con,x}(\tilde x)=V_j^{con,x}\cap(V_i^\perp+\{\tilde x\})$
over $\tilde x$. Transversality of $\tau_x$ with $V_i^\perp$
and uniqueness of $\lim T_{y_n}V_j^{con,x}$ imply that 
$V_j^{con,x}(\tilde x)$ is a manifold and $\tau_j(y)=T_yV_j\cap V_i^\perp$ 
depends continuously on $y$ in $V_j^{con,x}$. 
Consider the set of $\tilde x\in V_i$ for which have the corresponding
``level'' set $V_j^{con,x}(\tilde x)$ has $\tilde x$ in the closure, i.e. 
$\tilde x\in \overline{V_j^{con,x}(\tilde x)}$. Since $V_j^{con,x}$ is 
essential, the set of such $\tilde x$'s is everywhere dense
in a neighborhood of $x$ in $V_i$. Moreover, the ``angle'' function 
$\P^b$ is bounded in absolute value by $\dt$ on each local connected
``level'' component of $V_j^{con,x}(\tilde x)$ having $\tilde x$
in its closure. Thus, one can find a sequence of points 
$\{y_n\}\subset V_j^{con,x}$ tending to $x$ each point $y_n$ of which
belongs to a ``level'' connected component of 
$V_j^{con,x}(\pi^\perp_i(y_n))$, having $\pi^\perp_i(y_n)\in V_i$ in 
the closure. By construction $|\P^b(y_n)|<\dt$ for all $n$. Q.E.D.

\subsection{Separation of Planes}
Consider the real case. The complex case can be done 
in a similar way. Let $\tau_0$ and $\tau_1$ be two distinct 
orientable $k$-dimensional planes in $\R^m$. 
An orientable $(m-k)$-dimensional plane $l$ in $\R^m$ 
{\it separates} $\tau_0$ and $\tau_1$ if $l$ is transversal 
to $\tau_0$ and $\tau_1$ and the orientations induced by $\tau_0+l$ 
and $\tau_1+l$ in $\R^m$ are different. Notice that there always 
exists an open set of orientable $(m-k)$-planes separating 
any two distinct orientable $k$-plane.

{\bf Rolle's Lemma.} {\it If a continuous family of orientable 
$k$-planes $\{\tau_t \}_{t \in [0,1]}$ connects $\tau_0$ 
and $\tau_1$ and an orientable $(m-k)$-plane $l$ separates 
$\tau_0$ and $\tau_1$. Then for some $t^*\in (0,1)$ 
transversality of $\tau_{t^*}$ and $l$ fails.}

In what follows we use the transversality theorem \cite{GM} which 
says :\ {\it if $V \subset \R^m$ is a manifold, then almost every
plane of dimension $k$ is transversal to $V$}.  

\subsection{A reduction lemma.}\label{reduction} 
\blm \label{two} Let $V_j$ and $V_i$ be a distinct strata and 
$\dim V_j>\dim V_i$. Then there is a set of strata 
$\{V_j^p\}_{p\in \Z}$ (resp. $\{V_i^p\}_{p\in \Z}$) in $V_j$ 
(resp. in $V_i$) each of positive codimension in $V_j$ (resp. in
$V_i$) such that
\beq \label{include}
Sing_{a(\textup{resp. b})}(V_j,V_i) \subset \ \ 
\bigcup_{p\in \Z} Sing_{a(\textup{resp. b})}(V_j^p,V_i)
\bigcup_{p\in \Z} V_i^p.
\eneq
{\bf Remarks.} \ 
1. Inductive application of this lemma to the right-hand side of
(\ref{include}) reduces dimensions of $V_j^p$'s up to $\dim V_i$.

2. By the result of Lojasiewicz \cite{Lo} dimension of the
frontier of a semivariety ($Sing_{a(\textup{resp.} b)}(V_j^p,V_i)
\subset V_i\cap \overline{V^p_j}$)
has dimension strictly smaller that a semivariety itself.

3. By lemma \ref{one} the set $Sing_{a(\textup{resp. b})}(V_j,V_i)$ is 
a semivariety. Since a countable union of semivarieties of positive 
codimension in $V_i$ contains $Sing_{a(\textup{resp. b})}(V_j,V_i)$,  
$Sing_{a(\textup{resp. b})}(V_j,V_i)$ has a positive 
codimension in $V_i$ which proves Theorem \ref{sing}.
\elm

{\it Proof:} \ \ \ 
If $x \in Sing_{a(\textup{resp. b})}(V_j,V_i)$, then by the
construction of $\P^{a(\textup{resp.} b)}$, for some $\eps>0$ 
there is a sequence $\{y_n\}\subset V_j^{con,x}$ with 
$\P^{a(\textup{resp.} b)}(y_n)\to \eps$. There are two cases: 
                                     
$1$) there are different limits: 
$\P^{a(\textup{resp.} b)}(y'_n)\to \eps'$, 
$\P^{a(\textup{resp.} b)}(y_n)\to \eps''$, and $\eps'\neq \eps''$;
 
$2$) the limit $\P^{a(\textup{resp.} b)}(y_n)$ is unique, positive, 
and independent of $\{y_n\}$ \footnote{one can show that 
this case is impossible}.

Consider case $1$). By Sard's lemma there is a regular 
value $\eps^*\in (\eps',\eps'')$ of $\P^{a(\textup{resp.} b)}$. 
By the rank theorem 
$V_j^{\eps^*}=(\P^{a(\textup{resp.} b)})^{-1}(\eps^*)$ 
is a smooth semivariety of codimension $1$ in $V_j$. Let's show 
that $x \in \overline{V_j^{\eps^*}}$. Consider a local connected 
component $V_j^{con,x}$ and two sequences $\{y_n'\}$ and $\{y_n''\}$ 
in $V_j^{con,x}$ converging to $x$ such that 
$\P^{a(\textup{resp.} b)}(y_n')\to \eps'$ and 
$\P^{a(\textup{resp.} b)}(y_n'')\to \eps''$ as $n\to \infty$. 
$\P^{a(\textup{resp.} b)}$ is continuous and $V_j^{con,x}$ is
connected, thus we can connect each $y_n'$ and $y_n''$ in $V_j^{con,x}$ 
by a curve and find a sequence $\tilde y_n \to x$ for which 
$\P^{a(\textup{resp.} b)}(\tilde y_n)= \eps^*$. Thus 
$x\in \overline{V_j^{\eps^*}}$. Consider a countable dense set 
$\{\eps_p\}_{p\in \Z_+}$ in $[0,k+1]$ of regular values of 
$\P^{a(\textup{resp.} b)}$ so that for any two $\eps'\neq \eps''$, 
there is a separating $\eps_p \in (\eps',\eps'')$. Define 
$V^p_j=(\P^{a(\textup{resp.} b)})^{-1}(\eps_p)$. Thus any
$b$-irregular point $x$ is in the closure of the union 
$\cup_{p\in \Z_+} V^p_j$. After consideration of case $2$), we will 
prove that $V_j^p$ is $b$-irregular over $V_i$ at those $x$. 

Consider case $2$). By Lemma \ref{one} in this case if 
$x\in Sing_{a(\textup{resp. b})}(V_j,V_i)$, then $x$ belongs to $V_i\setminus
Un_a(V_j,V_i)$. Therefore, there are two sequences $\{y^0_n\},\ 
\{y^1_n\}$ in a local connected component $V_j^{con,x}$ tending to 
$x$ such that $T_{y^0_n}V_j \to \tau_0,\ T_{y^1_n}V_j \to \tau_1$, 
and $\tau_0 \neq \tau_1$. Choose an orientation of $T_{y^0_0}V_j$. 
By connecting $y^0_0$ locally with all other points $\{y^s_n\}$ one
can induce an orientation on all other $T_{y^s_n}V_j$ so that 
the orientations of $\tau_0$ and $\tau_1$ coincide with the
orientations of the limits.

Denote $\dim V_j$ by $k$. There is an orientable $(m-k)$-plane $l_j$ 
separating $\tau_0$ and $\tau_1$ and transversal to $V_j$ (by the
transversality theorem). Consider the orthogonal projection
$\pi_{l_j}$ along $l_j$ onto its orthogonal complement $l_j^\perp$.
Denote by $p_{l_j,j}$ its restriction to $V_j$, 
$p_{l_j,j}=\pi_{l_j}|_{V_j}:V_j \to l_j^\perp$.
Denote by $Crit(l_j,V_j)$ the set of critical points of 
$p_{l_j,j}$ in $V_j$ where the rank of $p_{l_j,j}$ is not maximal.
Then $Crit(l_j,V_j)$ is a semivariety in $V_j$ and 
$\dim Crit(l_j,V_j)<\dim V_j$. Connect two points $y_n^0$ and $y_n^1$ 
by a curve in $V_j$, then $T_{y_n^0}V_j$ deformates continuously  
to $T_{y_n^0}V_j$. Then  by Rolle's Lemma there is a critical point 
of $p_{l_j,j}$ in $V_j^{con,x}$ arbitrarily close to $x$. Thus
$x \in \overline{Crit(l_j,V_j)}$.

By the transversality theorem there is a countable dense set of 
orientable $(m-k)$-planes $\{l^r_j\}_{r \in \Z_+}$ transversal to 
$V_j$ and separating any two distinct orientable $k$-planes $\tau_0$ 
and $\tau_1$. Therefore, we have 
\beq \label{unique}
V_i \setminus Un_a(V_j,V_i) \subset \bigcup_{r\in \Z_+} 
\left\{\overline{Crit(l^r_j,V_j)}\setminus Crit(l^r_j,V_j)\right\}.
\eneq
By lemma \ref{one}  we know that $V_i \setminus Un_a(V_j,V_i)$ is 
a semivariety. We know that $Crit(l_j,V_j)$ $\subset V_j$ is a
semivariety and $\dim Crit(l_j,V_j) <\dim V_j$. Thus we can decompose
it into strata $Crit(l_j,V_j)=\bigsqcup_{p \in L_j} V^p_j$.
Renumerate these $V^p_j$'s to have $\{V^p_j\}_{p\in \Z_-}$.

Consider strata $\{V^p_j\}_{p\in \Z}\subset V_j$ which 
we constructed in the cases $1$ and $2$. Then  
\beq
Sing_{a(\textup{resp. b})}(V_j,V_i) \subset \bigcup_{p \in \Z}
\left\{\overline{V_j^p}\setminus V_j^p\right\}.
\eneq
The definitions of $\P^{a(\textup{resp.} b)}$ and $\pi_{j,s}$
explicitly imply that (\ref{include}) is satisfied, because 
$\P^{a(\textup{resp.} b)}(y_n)$ has a positive limit point for any 
$\{y_n\}\subset V_j^p$. If one projects along a smaller plane 
($T_{y_n}V^p_j\subset T_{y_n}V_j$), then the size of the projection 
is larger. Thus for the Kuo map 
$\P_{j^p,i}^{a(\textup{resp.} b)}:V^p_j\to \R$, defined in (\ref{kuofunc}),
the sequence $\P_{j^p,i}^{a(\textup{resp.} b)}(y_n)$ also has a
positive limit point. Now to separate interior and boundary points of 
the closures $\overline{V_j^p}$ in $V_i$ define the set
$V_i^p=(V_i\cap \overline{V_j^p})\setminus 
Int_{V_i}(V_i\cap \overline{V_j^p})$.
This completes the proof of the lemma and Theorem \ref{sing}. Q.E.D.

{\it Acknowledgments:}\ \ I would like to thank my thesis advisor 
John Mather and David Nadler for stimulating discussions and numerous 
remarks on mathematics and the English usage of the paper.


\begin{thebibliography}{Gned}
\def\bitem#1{\bibitem[#1]{#1}}
\bitem{GM} Goresky, M. MacPherson, R.; Stratified 
Morse Theory, Springer, 1987;

\bitem{Hi} Hironaka, H.; Number Theory, Algebraic
Geometry and Commutative Algebra, Volume in Honor of
Y. Akizuki, Kinokunia, Tokyo, 1973;

\bitem{Ja} Jakobson, N.; Basic Algebra, vol. 1, 1974; 

\bitem{Ku} Kuo, T.-C.; The ratio test for analytic Whitney
stratifications, Lecture Notes, No. 192, pp.141-149;

\bitem{Lo} Lojasiewicz, S.; Ensemble Semi-Analytiques,
IHES Lecture Notes, 1965;

\bitem{Ma} Mather, J.; Notes on Topological Stability,
Lecture Notes, Harvard University, 1970;

\bitem{Mi} Milnor, J.; Singularities of Complex 
Hypersurfaces, Ann. of Math. Studies, no. 61, 1968;

\bitem{Mu} Mumford, D.; Algebraic Geometry I, Spri\-n\-ger, 
New York, 1976

\bitem{PW} du Plessis, A. Wall. T.; The Geometry
of Topological Stability, Oxford, 1995;

\bitem{Th} Thom, R.; Propri\'et\'e Diff\'erentielle 
Locales des Ensembles Analytiques, Seminaire Bourbaki,
1964/65, exp. 281;

\bitem{Wa} Wall, T.; Regular Stratifications,
Lecture Notes in Mathematics, No. 468, pp. 332-344;

\bitem{Wh} Whitney, H.; Tangents to an Analytic Variety,
Ann. of Math. 81 (1965), pp. 496-549.
\end{thebibliography}
\end{document}